\documentclass[11pt]{article}
\usepackage{amsfonts}
\usepackage{amsmath}
\usepackage{amssymb}
\usepackage{amsthm}
\usepackage{xypic}
\usepackage[top=1in,bottom=1in,left=1.25in,right=1.25in]{geometry}

\date{}
\allowdisplaybreaks
\begin{document}

\renewcommand{\baselinestretch}{1.2}
\renewcommand{\arraystretch}{1.0}

\title{\bf Representation Crossed Category of Group-cograded Multiplier Hopf Algebras}
\author
{
  \textbf{Tao Yang} \footnote{Nanjing Agricultural University, Nanjing 210095, Jiangsu, CHINA.
             E-mail: tao.yang@njau.edu.cn}
}
\maketitle

\begin{center}
\begin{minipage}{12.cm}

 \textbf{Abstract}
 Let $A=\bigoplus_{p\in G}A_{p}$ be a multiplier Hopf $T$-coalgebra over a group $G$, 
 in this paper we give the definition of the crossed left $A$-$G$-modules
 and show that the category of crossed left $A$-$G$-modules is a monoidal category. 
 Finally we show that a family of multipliers $R = \{R_{p, q} \in M(A_{p}\otimes A_{q})\}_{p, q\in G}$ 
 is a quasitriangular structure of a multiplier $T$-coalgebra $A$ if and only if 
 the crossed left $A$-$G$-module category over $A$ is a braided monoidal category with the braiding $c$ defined by $R$,
 generalizing the main results in \cite{ZCL11} to the more general framework of multiplier Hopf algebras.

\\


\end{minipage}
\end{center}
\normalsize

\section{Introduction}
\def\theequation{\thesection.\arabic{equation}}
\setcounter{equation}{0}

 The motivating example for quasitriangular Hopf algebras is given by the Hopf algebra $U_{q}(\mathrm{g})$ 
 for $\mathrm{g}$ a finite-dimensional semisimple Lie algebra over the field $\mathbb{C}$ of complex numbers (see \cite{D87}).
 However, $H$ is not quasitriangular in the strict sense of the definition. 
 The $R$-matrix lies in a completion of $H\otimes H$ rather than in $H\otimes H$ itself. 
 The $U_{q}(\mathrm{g})$ is called 'topologically' quasitriangular.
 
 The explicit construction of the universal R-matrix is complicated. 
 There are some ways to evade this problem. One can study the finite-dimensional modules to obtain solutions of the Yang–Baxter equation.
 The present approach with multiplier Hopf algebras gives an alternative way to construct a generalized R-matrix in purely algebraic terms. 
 The notion of a quasitriangular (resp. group-cograded) multiplier Hopf algebra is introduced in \cite{Z99} (resp. \cite{DVW05}).

 In 2011, Zhu, Chen and Li \cite{ZCL11} studied the relation between the quasitriangular structure of a  Hopf group-coalgebra and the braided monoidal category.
 And in \cite{G15}, the author considered this relation in the weak Hopf group-coalgebra case.
 It is natural to ask whether we study the relation between the notion of quasitriangular group-cograded multiplier Hopf algebras 
 and its representation crossed categories.
 
 This question motivates the present research, and in this paper we will give a positive answer to the above question.
 
 The paper is organized in the following way.
 In section 2, we recall some notions which we will use in the following, such as the definition of
 multiplier Hopf algebra, multiplier Hopf $T$-coalgebras and their quasitriangular structures.
 In section 3, we introduce the notions of (crossed) left $A$-$G$-modules, 
 and show that the categories of (crossed) left $A$-$G$-modules are both monoidal categories.
 In section 4, we show the relation between multiplier Hopf $T$-coalgebra's quasitriangular structure and its representation crossed categories, 
 and prove that the monoidal category $({}_{A}\mathcal{M}_{crossed}, \otimes, K, a, l, r)$ is a braided monoidal category with the braiding $c$
 if and only if $A$ is a quasitriangular multiplier Hopf $T$-coalgebra, where $c$ is defined by $R$.

\section{Preliminaries}
\def\theequation{\thesection.\arabic{equation}}
\setcounter{equation}{0}

 Throughout this paper, all spaces we considered are over a fixed field $k$.
 Let $A$ be an (associative) algebra. We do not assume that $A$ has a unit, but we do require that
 the product, seen as a bilinear form, is non-degenerated. This means that, whenever $a\in A$ and $ab=0$ for all $b\in A$
 or $ba=0$ for all $b\in A$, we must have that $a=0$.
 Then we can consider the multiplier algebra $M(A)$ of $A$.
 Recall that $M(A)$ is characterized as the largest algebra with identity containing $A$ as an essential two-sided ideal.
 In particularly, we still have that, whenever $a\in M(A)$ and $ab=0$ for all $b\in A$ or $ba=0$ for all $b\in A$, again $a=0$.
 Furthermore, we consider the tensor algebra $A\otimes A$. It is still non-degenerated and we have its multiplier algebra $M(A\otimes A)$.
 There are natural imbeddings
 $$A\otimes A \subseteq M(A)\otimes M(A) \subseteq M(A\otimes A).$$
 In generally, when $A$ has no identity, these two inclusions are stict.
 If $A$ already has an identity, the product is obviously non-degenerate
 and $M(A)=A$ and $M(A\otimes A) = A\otimes A$. More details about the concept of the multiplier algebra of an algebra, we refer to \cite{V94}.

 Let $A$ and $B$ be non-degenerate algebras, if homomorphism $f: A\longrightarrow M(B)$ is non-degenerated
 (i.e., $f(A)B=B$ and $Bf(A)=B$),
 then has a unique extension to a homomorphism $M(A)\longrightarrow M(B)$, we also denote it $f$.

\subsection{Multiplier Hopf algebras}

 Now, we recall the definition of a multiplier Hopf algebra (see \cite{V94} for details).
 A comultiplication on algebra $A$ is a homomorphism $\Delta: A \longrightarrow M(A \otimes A)$ such that $\Delta(a)(1 \otimes b)$ and
 $(a \otimes 1)\Delta(b)$ belong to $A\otimes A$ for all $a, b \in A$. We require $\Delta$ to be coassociative in the sense that
 \begin{eqnarray*}
 (a\otimes 1\otimes 1)(\Delta \otimes \iota)(\Delta(b)(1\otimes c))
 = (\iota \otimes \Delta)((a \otimes 1)\Delta(b))(1\otimes 1\otimes c)
 \end{eqnarray*}
 for all $a, b, c \in A$ (where $\iota$ denotes the identity map).

 A pair $(A, \Delta)$ of an algebra $A$ with non-degenerate product and a comultiplication $\Delta$ on $A$ is called
 a \emph{multiplier Hopf algebra}, if the linear map $T_{1}, T_{2}$ defined by
 \begin{eqnarray}
 T_{1}(a\otimes b)=\Delta(a)(1 \otimes b), \qquad T_{2}(a\otimes b)=(a \otimes 1)\Delta(b)
 \end{eqnarray}
 are bijective.

 The bijectivity of the above two maps is equivalent to the existence of a counit and an antipode S satisfying (and defined by)
 \begin{eqnarray}
 && (\varepsilon\otimes\iota)(\Delta(a)(1\otimes b)) = ab, \qquad  m(S\otimes\iota)(\Delta(a)(1\otimes b))=\varepsilon(a)b, \label{1.1} \\
 && (\iota\otimes\varepsilon)((a\otimes 1)\Delta(b)) = ab, \qquad  m(\iota\otimes S)((a\otimes 1)\Delta(b))=\varepsilon(b)a,\label{1.2}
 \end{eqnarray}
 where $\varepsilon:A\longrightarrow k$ is a homomorphism, $S: A\longrightarrow M(A)$ is an anti-homomorphism
 and $m$ is the multiplication map, considered as a linear map from $A\otimes A$ to $A$ and extended to
 $M(A)\otimes A$ and $A\otimes M(A)$.

 A multiplier Hopf algebra $(A, \Delta)$ is called \emph{regular} if $(A, \Delta^{cop})$ is also a multiplier Hopf algebra,
 where $\Delta^{cop}$ denotes the co-opposite comultiplication defined as $\Delta^{cop}=\tau \circ \Delta$ with $\tau$ the usual flip map
 from $A\otimes A$ to itself (and extended to $M(A\otimes A)$). In this case, $\Delta(a)(b \otimes 1), (1 \otimes a)\Delta(b) \in A \otimes A$
 for all $a, b\in A$.
 By Proposition 2.9 in \cite{V98}, multiplier Hopf algebra $(A, \Delta)$ is regular if and only if the antipode $S$ is bijective from $A$ to $A$.

 Remark that any Hopf algebra is a multiplier Hopf algebra and conversely, any multiplier Hopf algebra, with identity, is a Hopf algebra.
 We will use the adapted Sweedler notation for multiplier Hopf algebras (see \cite{V08}), e.g., write $a_{(1)} \otimes a_{(2)}b$ for $\Delta(a)(1 \otimes b)$ 
 and $ab_{(1)} \otimes b_{(2)}$ for $(a \otimes 1)\Delta(b)$.

 Let $A$ be a regular multiplier Hopf algebra. Suppose $X$ is a left $A$-module
 with the module structure map $\cdot : A\otimes X\longrightarrow X$. 
 We will always assume that the module is non-degenerate, this means that $x=0$ if $x\in X$ and $a\cdot x=0$ for all $a\in A$.
 If the module is unital (i.e., $A\cdot X=X$), then we can get an extension of the module structure to $M(A)$, 
 this means that we can define $f\cdot x$, where $f\in M(A)$ and $x\in X$.
 In fact, since $x\in X=A\cdot X$, then $x= \sum_{i} a_{i}\cdot x_{i}$, $f\cdot v=\sum_{i} (fa_{i})\cdot x_{i}$.
 In this setting, we can easily get $1_{M(A)}\cdot x=x$.

 \subsection{Multiplier Hopf $T$-coalgebras and their quasitriangular structures}

 Let $(A, \Delta)$ be a multiplier Hopf algebra and $G$ a group with unit $e$.
 Assume that there is a family of (non-trivial) subalgebras $\{A_{p}\}_{p\in G}$ of $A$ so that
 \begin{enumerate}
 \item[(1)] $A=\bigoplus_{p\in G}A_{p}$ with $A_{p}A_{q}=0$ whenever $p,q\in G$ and $p\neq q$.
 \item[(2)] $\Delta(A_{pq})(1\otimes A_{q})=A_{p}\otimes A_{q}$ and $(A_{p}\otimes 1)\Delta(A_{pq})=A_{p}\otimes A_{q}$ for all $p,q \in G$.
 \end{enumerate}
 Then $(A, \Delta)$ is called  a \emph{$G$-cograded multiplier Hopf algebra}(see \cite{DVW05,YW11a}).

 We extend the Sweedler notation for a comultiplication in the following way:
 for any $p, q \in G$, $a\in A_{pq}$ and $a'\in A_{q}$, we write
 \begin{eqnarray*}
 \Delta_{p, q}(a)(1\otimes a') = a_{(1, p)} \otimes a_{(2, q)}a'.
 \end{eqnarray*}

 Let $A$ be a $G$-cograded multiplier Hopf algebra, then $A$ has the form $A=\bigoplus_{p\in G} A_{p}$.
 Assume that there is a group homomorphism $\pi: G\longrightarrow Aut(A)$.
 We call $\pi$ an \emph{admissible} action of $G$ on $A$ if also the following
 requirements hold
 \begin{enumerate}
 \item[(1)] $\Delta(\pi_{p}(a))=(\pi_{p}\otimes \pi_{p})\Delta(a)$ for all $a\in A$.
 \item[(2)] $\pi_{p}(A_{q})=A_{\rho_{p}(q)}$, where $\rho$ is an action of the group $G$ on itself.
 \item[(3)] $\pi_{\rho_{p}(q)}=\pi_{pqp^{-1}}$. \\
            This means that the map $\pi$ takes care of $\rho$ not being the adjoint action.
            If $\rho$ is the adjoint action, $\pi$ is called a \emph{crossing}.
 \end{enumerate}

 A group-cograded multiplier Hopf algebra $A=\bigoplus_{p\in G} A_{p}$ is said to be a \emph{multiplier Hopf $T$-coalgebra}
 provided it is endowed with a crossing $\pi$ such that each $\pi_{q}$ preserves the comultiplication and the counit,
 i.e., for all $p, q, r \in G$,
 \begin{eqnarray*}
 (\pi_{q} \otimes \pi_{q}) \Delta_{p, r} = \Delta_{qpq^{-1}, qrq^{-1}} \pi_{q}, \quad
 \varepsilon \pi_{q} = \varepsilon,
 \end{eqnarray*}
 and $\pi$ is multiplicative in the sense that $\pi_{pq}=\pi_{p} \pi_{q}$ for all $p, q \in G$.
 It can be considered as generalization of crossed Hopf $G$-coalgebra introduced in \cite{T08}.
 
 Let $A$ be a multiplier Hopf $T$-coalgebra, then we can construct a new regular multiplier Hopf algebra on $A$ 
 by deforming the comultiplication while the algebra structure on $A$ is kept (see Theorem 3.11 in \cite{DV07}).
 The comultiplication deformation of $A$ depends on the crossing $\pi$ in the following way:
 for all $a\in A$ and $a'\in A_{q}$,
 \begin{eqnarray*}
 \widetilde{\Delta}(a)(1\otimes a')=(\pi_{q^{-1}}\otimes \iota)(\Delta(a)(1\otimes a')).
 \end{eqnarray*}
 
 Recall from \cite{DVW05}, a $G$-cograded multiplier Hopf algebra with a crossing action $\pi$ is called {\it quasitriangular}
 if there is a multiplier $R=\sum _{p, q \in G}R_{p, q}$
 with  $R_{p, q} \in M(A_{p}\otimes A_{q})$  such that
 \begin{eqnarray*}
 (\pi_{p}\otimes \pi_{p})(R) = R,  && R\Delta (a)=(\widetilde{\Delta} )^{cop}(a)R, \\
 (\widetilde{\Delta}\otimes  \iota)(R)=R_{13}R_{23}, && (\iota \otimes \Delta)(R) = R_{13} R_{12}.
 \end{eqnarray*}
 for all $p\in G$ and $a\in A$. Sometimes we call $R$ a generalized $R$-matrix.

\section{Category of Crossed Left Modules}
\def\theequation{\thesection.\arabic{equation}}
\setcounter{equation}{0}

 The concept of a group-cograded multiplier Hopf algebra was introduced by Abd El-hafez et al in \cite{ADV07}
 as a generalization of Hopf group-coalgebras introduced in \cite{T08}. 
 In the following, similar to the Hopf group-coalgebra case, 
 we will introduce the crossed module category for a $G$-cograded multiplier Hopf algebra $A=\bigoplus_{p\in G} A_{p}$.
 \\
 
 \textbf{Definition \thesection.1}
 Let $A$ be a $G$-cograded multiplier Hopf algebra, a left $A$-$G$-module is a family $M = \{M_{p}\}_{p\in G}$ of vector spaces such that 
 $M_{p}$ is a left unital $A_{p}$-module for any $p\in G$. 
 We denote the structure maps of left $A_{p}$-module $M_{p}$ and $A$-$G$-module $M$ by $\Gamma_{M_{p}}: A_{p}\otimes M_{p} \rightarrow M_{p}$
 and $\Gamma_{M} = \{\Gamma_{M_{p}}\}_{p\in G}$, respectively.
 Furthermore, if $\{M_{p}\}_{p\in G}$ and $\{N_{p}\}_{p\in G}$ are two left $A$-$G$-modules, 
 a left $A$-$G$-morphism is a family $f=\{f_{p}: M_{p} \rightarrow N_{p}\}_{p\in G}$
 of linear maps such that $f_{p}$ is an left $A_{p}$-morphism for any $p\in G$.
 \\
 
 Denote by ${}_{A}\mathcal{M}$ the category of all left $A$-$G$-modules and its morphisms are left $A$-$G$-morphisms.
 Let $\{M_{p}\}_{p\in G}$ and $\{N_{p}\}_{p\in G}$ be two left $A$-$G$-modules, we can easily check that 
 $M_{s}\otimes N_{t}$ is a unital $A_{s}\otimes A_{t}$-module, and a unital left $A_{st}$-module with the action 
 \begin{eqnarray*}
 \Gamma_{M_{s}\otimes N_{t}}(a\otimes (m\otimes n)):= a\cdot (m\otimes n) = a_{(1, s)}\cdot m \otimes a_{(2, t)}\cdot n,
 \end{eqnarray*}
 for all $a\in A_{st}$, $m\in M_{s}$ and $n\in N_{t}$. 
 Hence, $(M\otimes N)_{p} := \bigoplus_{st=p} M_{s}\otimes N_{t}$ is a unital left $A_{p}$-module.
 Thus $M\otimes N = \{(M\otimes N)_{p}\}_{p\in G}$ is a left $A$-$G$-module with the structure map 
 $\Gamma_{M\otimes N} = \{\Gamma_{(M\otimes N)_{p}}\}_{p\in G}$ given by
 \begin{eqnarray}
 \Gamma_{(M\otimes N)_{p}} 
 = \bigoplus_{st=p} (\Gamma_{M_{s}} \otimes \Gamma_{N_{t}}) (\iota\otimes\tau\otimes\iota_{N_{t}}) 
   (\Delta_{s, t}\otimes\iota_{M_{s}}\otimes\iota_{N_{t}}). \label{1}
 \end{eqnarray}
 
 Suppose that $L=\{L_{p}\}_{p\in G}$ is also a left $A$-$G$-module, then we have two left $A$-$G$-modules $(L\otimes M)\otimes N$ and $L\otimes (M\otimes N)$.
 By the $A$-$G$-module action (\ref{1}), we have for any $p\in G$
 \begin{eqnarray*}
 ((L\otimes M)\otimes N)_{p} 
 &=& \bigoplus_{st=p} (L\otimes M)_{s}\otimes N_{t} 
 = \bigoplus_{st=p} (\bigoplus_{qr=s} L_{q}\otimes M_{r}) \otimes N_{t} \\
 &=& \bigoplus_{qrt=p}  (L_{q}\otimes M_{r}) \otimes N_{t}, \\
 (L\otimes (M\otimes N))_{p}
 &=& \bigoplus_{qs=p} L_{q}\otimes (M\otimes N)_{s} 
 = \bigoplus_{qs=p} L_{q}\otimes (\bigoplus_{rt=s} M_{r}\otimes N_{t}) \\
 &=& \bigoplus_{qrt=p}  L_{q}\otimes (M_{r} \otimes N_{t}).
 \end{eqnarray*}
 One can easily check that 
 $a_{q, r, t}: (L_{q}\otimes M_{r}) \otimes N_{t} \rightarrow L_{q}\otimes (M_{r} \otimes N_{t}), 
 (l\otimes m)\otimes n \mapsto l\otimes (m\otimes n)$ is an isomorphism of unital $A_{qrt}$-modules.
 Then for any $p\in G$, $a_{p} = \bigoplus_{qrt=p} a_{q, r, t}$ is an isomorphism of left $A_{p}$-modules 
 from $((L\otimes M)\otimes N)_{p}$ to $(L\otimes (M\otimes N))_{p}$, and 
 $a=\{a_{p}\}_{p\in G}: (L\otimes M)\otimes N \rightarrow L\otimes (M\otimes N)$ is a left $A$-$G$-module isomorphism,
 it is a family of natural isomorphisms.

 Since $A_e$ is a multiplier Hopf algebra, the basic field $k$ is a unital left $A_e$-module with the action $a\cdot \kappa = \varepsilon(a)\kappa$,
 where $a\in A_e$ and $\kappa \in k$.
 Let $K=\{K_p\}_{p\in G}$ with $K_e=k$ and $K_p=0$ whenever $p\neq e$, then $K$ is a left $A$-$G$-module.
 For any left $A$-$G$-module $M$, $(K\otimes M)_p = K_e \otimes M_{p} = k\otimes M_{p}$ and 
 $(M\otimes K)_p = M_{p}\otimes K_e = M_p\otimes k$. So we have natural isomorphisms $l_M: K\otimes M\rightarrow M$ and $r_M: M\otimes K\rightarrow M$
 defined by 
 \begin{eqnarray*}
 (l_M)_p: k\otimes M_{p}\rightarrow M_{p}, && \kappa\otimes m\mapsto m, \\
 (r_M)_p: M_{p}\otimes k\rightarrow M_{p}, && m\otimes\kappa\mapsto m.
 \end{eqnarray*} 
 That is, $\{l_M\}$ and $\{r_M\}$ are two families of natural isomorphisms of left $A$-$G$-modules.
 
 From the above discussion, We summarize the following result:
 \\
 
 \textbf{Theorem \thesection.2}
 $({}_{A}\mathcal{M}, \otimes, K, a, l, r)$ is a monoidal category with the unit object $K$.
 \\
 
 In the following, we will consider the crossed module of multiplier Hopf $T$-coalgebra $A=\bigoplus_{p\in G} A_{p}$. 
 Firstly, we introduce the definition of crossed left $A$-$G$-module, which generalizes the notion in \cite{ZCL11}.
  
 \textbf{Definition \thesection.3}
 Let $A$ be a multiplier Hopf $T$-coalgebra, a left $A$-$G$-module $M = \{M_{p}\}_{p\in G}$ is called crossed if there exists 
 a family $\pi_{M} = \{\pi_{M, p}: M_{q}\rightarrow M_{pqp^{-1}}\}_{p\in G}$ of linear maps such that the following conditions are satisfied
 \begin{enumerate}
 \item[(1)] each $\pi_{M, p}: M_{q}\rightarrow M_{pqp^{-1}}$ is a vector space isomorphism,
 \item[(2)] each $\pi_{M, p}$ preserves the action, ie., for any $p, q\in G$, $\pi_{M, p}(a\cdot m) = \pi_{p}(a)\cdot \pi_{M, p}(m)$,
 \item[(3)] $\pi_{M}$ is multiplicative in the sense that for any $p, q\in G$, $\pi_{M, p}\circ\pi_{M, q} = \pi_{M, pq}$.
 \end{enumerate}
 And for two left crossed $A$-$G$-modules $M=\{M_{p}\}_{p\in G}$ and $N=\{N_{p}\}_{p\in G}$, 
 a left $A$-$G$-module morphism $f=\{f_{p}\}_{p\in G}: M\rightarrow N$ is called crossed 
 if furthermore for any $p, q\in G$, $\pi_{N, q}f_{p} = f_{qpq^{-1}} \pi_{M, q}$.
 \\
 
 Let $A=\bigoplus_{p\in G} A_{p}$ be a multiplier Hopf $T$-coalgebra.
 Denote by ${}_{A}\mathcal{M}_{crossed}$ the category of all crossed left $A$-$G$-modules and crossed left $A$-$G$-morphisms.
 Let $M$ and $N$ be two crossed left $A$-$G$-modules, we have already shown that $M\otimes N$ is also a left $A$-$G$-module.
 Define $\pi_{M\otimes N, q}: (M\otimes N)_{p}\rightarrow (M\otimes N)_{qpq^{-1}}$ by
 \begin{eqnarray}
 \pi_{M\otimes N, q}|_{M_{s}\otimes N_{t}} = \pi_{M, q}|_{M_{s}} \otimes \pi_{N, q}|_{N_{t}}, \label{3}
 \end{eqnarray} 
 where $p, q, s, t\in G$ with $p=st$. Since
 \begin{eqnarray*}
 (M\otimes N)_{p} 
 &=& \bigoplus_{st=p} M_{s}\otimes N_{t}, \quad \mbox{and} \\
 (M\otimes N)_{qpq^{-1}} 
 &=& \bigoplus_{qstq^{-1}=qpq^{-1}} M_{qsq^{-1}}\otimes N_{qtq^{-1}} \\
 &=& \bigoplus_{st=p} M_{qsq^{-1}}\otimes N_{qtq^{-1}},
 \end{eqnarray*} 
 $\pi_{M\otimes N, q}$ is well-defined $k$-linear isomorphism from $(M\otimes N)_{p}$ to $(M\otimes N)_{qpq^{-1}}$.
 Moreover, for $a\in A_{p}$, $m\in M_{s}$ and $n\in N_{t}$ we have 
 $\pi_{M\otimes N, q}(a\cdot (m\otimes n)) = \pi_{q}(a)\cdot \pi_{M\otimes N, q}(m\otimes n) $.
 Indeed,
 \begin{eqnarray*}
 \pi_{M\otimes N, q}(a\cdot (m\otimes n)) 
 &=& (\pi_{M, q} \otimes \pi_{N, q})(a_{(1, s)} \cdot m\otimes a_{(2, t)} \cdot n) \\
 &=& \pi_{M, q}(a_{(1, s)} \cdot m)\otimes \pi_{N, q}(a_{(2, t)} \cdot n) \\
 &=& \pi_{q}(a_{(1, s)}) \cdot \pi_{M, q}(m)\otimes \pi_{q}(a_{(2, t)}) \cdot \pi_{N, q}(n),
 \end{eqnarray*} 
 and 
 \begin{eqnarray*}
 \pi_{q}(a)\cdot \pi_{M\otimes N, q}(m\otimes n)
 &=& \pi_{q}(a)\cdot (\pi_{M, q}(m)\otimes \pi_{N, q}(n)) \\
 &=& \Delta_{qsq^{-1}, qtq^{-1}}(\pi_{q}(a))\cdot (\pi_{M, q}(m)\otimes \pi_{N, q}(n)) \\
 &=& (\pi_{q}\otimes \pi_{q})\Delta_{s, t}(a)\cdot (\pi_{M, q}(m)\otimes \pi_{N, q}(n)) \\
 &=& \pi_{q}(a_{(1, s)}) \cdot \pi_{M, q}(m)\otimes \pi_{q}(a_{(2, t)}) \cdot \pi_{N, q}(n).
 \end{eqnarray*} 
 It is easy to get that $\pi_{M\otimes N, q}\circ\pi_{M\otimes N, q'} = \pi_{M\otimes N, qq'}$, 
 i.e., $\pi_{M\otimes N, q}$ is multiplicative.
 Thus $M\otimes N$ is a crossed left $A$-$G$-module.
 
 Now let $L$ be another crossed left $A$-$G$-module, one can easily check that 
 \begin{eqnarray*}
 \pi_{L\otimes(M\otimes N), q} a_{p} = a_{qpq^{-1}} \pi_{(L\otimes M)\otimes N, q},
 \end{eqnarray*} 
 and so $a: (L\otimes M)\otimes N\rightarrow L\otimes(M\otimes N)$ is a crossed left $A$-$G$-module morphism.

 Note that $K_e = k = K_{qeq^{-1}}$ and $K_p = 0 = K_{qpq^{-1}}$ for $p\neq e$. 
 Since $\pi_{e}: A_{p}\rightarrow A_{p}$ is the identity map, set $\pi_{K, q}= \iota: K_{p}\rightarrow K_{qpq^{-1}}$, 
 then the unit object $K=\{K_{p}\}_{p\in G}$ is a crossed left $A$-$G$-module.
 
 Finally, one can easily check that the left and right unit constraints $l=\{l_M\}$ and $r=\{r_M\}$ are crossed left $A$-$G$-module morphisms.
 
 Therefore we get another main result of this section:
 \\
 
 \textbf{Theorem \thesection.4}
 $({}_{A}\mathcal{M}_{crossed}, \otimes, K, a, l, r)$ is a monoidal category with the unit object $K$.

\section{Braided Monnidal Category}
\def\theequation{\thesection.\arabic{equation}}
\setcounter{equation}{0}

 Throughout the following, Assume that $A=\bigoplus_{p\in G} A_{p}$ be a multiplier Hopf $T$-coalgebra and 
 $R = \{R_{p, q} \in M(A_{p}\otimes A_{q})\}_{p, q\in G}$ be a family of multipliers.
 Let $M$ and $N$ be any two crossed left $A$-$G$-modules, for $s, t\in G$ we define
 \begin{eqnarray*}
 c_{M_{s}, N_{t}}: M_{s}\otimes N_{t} \longrightarrow N_{sts^{-1}}\otimes M_{s}
 \end{eqnarray*} 
 by
 \begin{eqnarray}
 c_{M_{s}, N_{t}}(m\otimes n) = (\pi_{N, s} \otimes \iota_{M_{s}}) \tau_{s, t} \Big(R_{s, t}\cdot (m\otimes n)\Big), \label{2}
 \end{eqnarray} 
 where $m\in M_{s}$ and $n\in N_{t}$. For any $p\in G$, define
 \begin{eqnarray*}
 (c_{M, N})_{p}: (M\otimes N)_{p} = \bigoplus_{st=p} M_{s}\otimes N_{t} \longrightarrow (N\otimes M)_{p} = \bigoplus_{st=p} N_{sts^{-1}}\otimes M_{s}
 \end{eqnarray*} 
 by 
 \begin{eqnarray*}
 (c_{M, N})_{p} = \bigoplus_{st=p} c_{M_{s}, N_{t}}.
 \end{eqnarray*} 
 Then it is obvious that $(c_{M, N})_{p}$ is a $k$-linear isomorphism for any $p\in G$ if and only if so is $c_{M_{s}, N_{t}}$ for any $s, t\in G$ with $st=p$.
 \\
 
 \textbf{Lemma \thesection.1} With the notations above, we have
 \begin{enumerate}
 \item[(1)] For any $p\in G$, $(c_{M, N})_{p}$ is a $k$-linear isomorphism for any crossed left $A$-$G$-modules $M$ and $N$ 
             if and only if $R$ is a family of invertible multipliers.
 \item[(2)] For any crossed left $A$-$G$-modules $M$ and $N$, $c_{M, N}: M\otimes N \rightarrow N\otimes M$ is a left $A$-$G$-module morphism if and only if
 \begin{eqnarray}
 R_{s, t}\Delta_{s, t}(a) = \Big(\tau_{t, s} (\pi_{s^{-1}}\otimes \iota)\Delta_{sts^{-1}, s}(a) \Big) R_{s, t}
 \end{eqnarray} 
 for $s, t\in G$ and $a\in A_{st}$.
 \end{enumerate}
 
 \emph{Proof}
 (1) If $R = \{R_{p, q} \in M(A_{p}\otimes A_{q})\}_{p, q\in G}$ is a family of invertible multipliers, 
 then obviously $c_{M_{s}, N_{t}}$ is a $k$-linear isomorphism for any $s, t\in G$ since so is $\pi_{N, s}$.
 Conversely, let $M = N = A$, because $\pi_{p}$ is an isomorphism, from the hypothesis one knows that the map
 $A_{s}\otimes A_{t}\rightarrow A_{s}\otimes A_{t}, x\otimes y\mapsto R_{s, t}(x\otimes y)$ is a $k$-linear isomorphism for any $s, t\in G$.
 It follows that $R$ is a invertible multiplier in the multiplier algebra $M(A_{s}\otimes A_{t})$.
 
 (2) It is clear that $(c_{M, N})_{p}$ is an $A_{p}$-module morphism if and only if 
 $c_{M_{s}, N_{t}}$ is an $A_{st}$-module morphism for any $s, t\in G$ with $st=p$.
 Let $m\in M_{s}, n\in N_{t}$ and $a\in A_{st}$, then we have
 \begin{eqnarray*}
 c_{M_{s}, N_{t}}(a\cdot (m\otimes n)) 
 &=& c_{M_{s}, N_{t}}(\Delta_{s, t}(a)\cdot (m\otimes n)) \\
 &\stackrel{(\ref{2})}{=}&  (\pi_{N, s} \otimes \iota_{M_{s}}) \tau_{s, t} \Big(R_{s, t}\cdot (\Delta_{s, t}(a)\cdot (m\otimes n))\Big) \\
 &=&  (\pi_{N, s} \otimes \iota_{M_{s}}) \tau_{s, t} \Big(R_{s, t}\Delta_{s, t}(a)\cdot (m\otimes n)\Big),
 \end{eqnarray*} 
 and
 \begin{eqnarray*}
 a\cdot c_{M_{s}, N_{t}}(m\otimes n) 
 &=& \Delta_{sts^{-1}, s}(a)\cdot c_{M_{s}, N_{t}}(m\otimes n) \\
 &\stackrel{(\ref{2})}{=}& (\pi_{s} \otimes \iota)\Big( (\pi_{s^{-1}} \otimes \iota)\Delta_{sts^{-1}, s}(a) \Big) 
     \cdot (\pi_{N, s} \otimes \iota_{M_{s}}) \tau_{s, t} \Big(R_{s, t}\cdot (m\otimes n) \Big) \\
 &=&  (\pi_{N, s} \otimes \iota_{M_{s}}) \Bigg((\pi_{s^{-1}} \otimes \iota)\Delta_{sts^{-1}, s}(a) 
     \cdot \tau_{s, t} \Big(R_{s, t} \cdot (m\otimes n)\Big) \Bigg) \\
 &=& (\pi_{N, s} \otimes \iota_{M_{s}}) \tau_{s, t} 
    \Bigg( \tau_{t, s}(\pi_{s^{-1}} \otimes \iota)\Delta_{sts^{-1}, s}(a) \cdot \Big(R_{s, t} \cdot (m\otimes n)\Big) \Bigg) \\
 &=& (\pi_{N, s} \otimes \iota_{M_{s}}) \tau_{s, t} 
    \Bigg( \Big(\tau_{t, s}(\pi_{s^{-1}} \otimes \iota)\Delta_{sts^{-1}, s}(a)\Big) R_{s, t} \cdot (m\otimes n) \Bigg).
 \end{eqnarray*} 
 Because $\pi_{N, s}$ is an isomorphism, $c_{M_{s}, N_{t}}(a\cdot (m\otimes n)) = a\cdot c_{M_{s}, N_{t}}(m\otimes n)$
 if and only if $R_{s, t}\Delta_{s, t}(a)\cdot (m\otimes n) =  \Big(\tau_{t, s}(\pi_{s^{-1}} \otimes \iota)\Delta_{sts^{-1}, s}(a)\Big) R_{s, t} \cdot (m\otimes n)$.
 Hence if $R_{s, t}\Delta_{s, t}(a) = \Big(\tau (\pi_{s^{-1}}\otimes \iota)\Delta_{sts^{-1}, s}(a) \Big) R_{s, t}$,
 then $c_{M_{s}, N_{t}}$ is an $A_{st}$-module morphism.
 Conversely, let $M = N = A$, then for any $x\in A_{s}$ and $y\in A_{t}$,
 $R_{s, t}\Delta_{s, t}(a)(x\otimes y) =  \Big(\tau_{t, s}(\pi_{s^{-1}} \otimes \iota)\Delta_{sts^{-1}, s}(a)\Big) R_{s, t} (x\otimes y)$.
 By the non-degenerate of the product, we can get
 $R_{s, t}\Delta_{s, t}(a) =  \Big(\tau_{t, s}(\pi_{s^{-1}} \otimes \iota)\Delta_{sts^{-1}, s}(a)\Big) R_{s, t}$.
 $\hfill \Box$
 \\
 
 \textbf{Lemma \thesection.2} The following two statements are equivalent:
 \begin{enumerate}
 \item[(1)] $\pi_{N\otimes M, q}\circ (c_{M, N})_{p} = (c_{M, N})_{qpq^{-1}} \circ \pi_{M\otimes N, q}$ 
 for any crossed left $A$-$G$-modules $M$ and $N$, and $p, q\in G$. 
 \item[(2)] $(\pi_{q}\otimes \pi_{q})R_{s, t} = R_{qsq^{-1}, qtq^{-1}}$ for any $q, s, t\in G$.
 \end{enumerate}
 
 \emph{Proof}
 For any $s, t, q\in G$, $m\in M_{s}$ and $n\in N_{t}$,
 \begin{eqnarray*}
 && \pi_{N\otimes M, q}\circ (c_{M, N})_{st}(m\otimes n) \\
 &=& (\pi_{N, q}\otimes \pi_{M, q}) (c_{M_{s}, N_{t}}) (m\otimes n) \\
 &\stackrel{(\ref{2})}{=}& (\pi_{N, q}\otimes \pi_{M, q}) (\pi_{N, s} \otimes \iota_{M_{s}}) \tau_{s, t} \Big(R_{s, t}\cdot (m\otimes n)\Big) \\
 &=& (\pi_{N, qs} \otimes \pi_{M, q})  \Big(\tau_{s, t}(R_{s, t})\cdot (n\otimes m)\Big) \\
 &=& (\pi_{N, qsq^{-1}} \otimes \iota_{M_{qsq^{-1}}}) (\pi_{N, q} \otimes \pi_{M, q}) \Big(\tau_{s, t}(R_{s, t})\cdot (n\otimes m)\Big) \\
 &=& (\pi_{N, qsq^{-1}} \otimes \iota_{M_{qsq^{-1}}}) \Big( (\pi_{q} \otimes \pi_{q}) (\tau_{s, t}R_{s, t})\cdot (\pi_{N, q}(n) \otimes \pi_{M, q}(m)) \Big) \\
 &=& (\pi_{N, qsq^{-1}} \otimes \iota_{M_{qsq^{-1}}}) 
    \Bigg( \tau_{qsq^{-1}, qtq^{-1}} (\pi_{q} \otimes \pi_{q}) (R_{s, t})\cdot \Big(\pi_{N, q}(n) \otimes \pi_{M, q}(m) \Big) \Bigg),
 \end{eqnarray*} 
 and 
 \begin{eqnarray*}
 && (c_{M, N})_{qstq^{-1}} \circ \pi_{M\otimes N, q}(m\otimes n) \\
 &=& c_{M_{qsq^{-1}}, N_{qtq^{-1}}} \big(\pi_{M, q}(m) \otimes \pi_{N, q}(n) \big) \\
 &\stackrel{(\ref{2})}{=}& 
      (\pi_{N, qsq^{-1}} \otimes \iota_{M_{qsq^{-1}}}) \tau_{qsq^{-1}, qtq^{-1}} 
      \Big(R_{qsq^{-1}, qtq^{-1}}\cdot \big(\pi_{M, q}(m) \otimes \pi_{N, q}(n) \big) \Big) \\
 &=& (\pi_{N, qsq^{-1}} \otimes \iota_{M_{qsq^{-1}}}) 
      \Big(\tau_{qsq^{-1}, qtq^{-1}} (R_{qsq^{-1}, qtq^{-1}}) \cdot \big(\pi_{N, q}(n) \otimes \pi_{M, q}(m) \big) \Big).
 \end{eqnarray*} 
 It follows that (2) implies (1). Conversely, assume that (1) is satisfied. 
 Because $\pi_{N, qsq^{-1}}$ is an isomorphism, it follows
 $(\pi_{q} \otimes \pi_{q}) (R_{s, t})\cdot (\pi_{M, q}(m) \otimes \pi_{N, q}(n))
 = R_{qsq^{-1}, qtq^{-1}} \cdot (\pi_{M, q}(m) \otimes \pi_{N, q}(n))$.
 Because $\pi_{M, q}$ and $\pi_{N, q}$ are isomorphisms, 
  $(\pi_{q} \otimes \pi_{q}) (R_{s, t})\cdot (u \otimes v) = R_{qsq^{-1}, qtq^{-1}} \cdot (u \otimes v)$ 
  for any $u\in M_{qsq^{-1}}$ and $v\in N_{qtq^{-1}}$.
 Let $M = N = A$, then by the non-degenerate of the product, we can get 
 $(\pi_{q} \otimes \pi_{q}) (R_{s, t}) = R_{qsq^{-1}, qtq^{-1}}$.
 $\hfill \Box$
 \\
 
 \textbf{Lemma \thesection.3}
  The following two statements hold:
 \begin{enumerate}
 \item[(1)] $c_{L, M\otimes N} = (\iota_{M}\otimes c_{L, N})(c_{L, M}\otimes \iota_{N})$ for any crossed left $A$-$G$-modules $L$, $M$ and $N$, 
 if and only if for any $r, s, t\in G$,
 \begin{eqnarray*}
 (\iota \otimes \Delta_{s, t})(R_{r, st}) = (R_{r, t})_{1s3} (R_{r, s})_{12t}.
 \end{eqnarray*}
 \item[(2)] $c_{L\otimes M, N} = (c_{L, N}\otimes \iota_{M})(\iota_{L}\otimes c_{M, N})$ for any crossed left $A$-$G$-modules $L$, $M$ and $N$,
 if and only if for any $r, s, t\in G$,
 \begin{eqnarray*}
 (\Delta_{r, s} \otimes \iota)(R_{rs, t}) = \Big((\iota\otimes \pi_{s^{-1}})R_{r, sts^{-1}}\Big)_{1s3} (R_{s, t})_{r23}.
 \end{eqnarray*} 
 \end{enumerate}
 
 \emph{Proof} 
 We only check the first statement and the second one is similar.
 For $v\in L_{r}$, $m\in M_{s}$ and $n\in N_{t}$,
 \begin{eqnarray*}
  && (c_{L, M\otimes N})_{rst}(v\otimes m\otimes n) \\
  &=& (\pi_{M\otimes N, r}\otimes \iota_{L_{r}}) \tau_{r, st} \Big(R_{r, st}\cdot (v\otimes (m\otimes n))\Big) \\
  &=& (\pi_{M\otimes N, r}\otimes \iota_{L_{r}}) (\iota_{M_{s}}\otimes \tau_{r, t}) (\tau_{r, s}\otimes \iota_{N_{t}}) 
     \Big((\iota\otimes\Delta_{s, t})R_{r, st}\cdot (v\otimes m\otimes n)\Big),
 \end{eqnarray*}
 and
 \begin{eqnarray*}
  && (\iota_{M}\otimes c_{L, N})(c_{L, M}\otimes \iota_{N})(v\otimes m\otimes n) \\
  &\stackrel{(\ref{2})}{=}&  (\iota_{M}\otimes c_{L, N}) \Big( (\pi_{M, r} \otimes \iota_{L_{r}}) \tau_{r, s} (R_{r, s}\cdot (v\otimes m)) \otimes n\Big) \\
  &=& (\pi_{M, r}\otimes c_{L, N}) (\tau_{r, s} \otimes \iota_{N_t}) \Big( (R_{r, s})_{12t}\cdot (v\otimes m \otimes n) \Big) \\
  &=& (\pi_{M, r}\otimes c_{L, N}) (R^{2}_{r, s}\cdot m\otimes R^{1}_{r, s}\cdot v \otimes n) \\
  &\stackrel{(\ref{2})}{=}& (\pi_{M, r}\otimes (\pi_{N, r}\otimes \iota_{L_{r}})) (\iota_{M_{s}}\otimes \tau_{r, t}) 
              \Big(R^{2}_{r, s}\cdot m\otimes R^{1}_{r, t}R^{1}_{r, s}\cdot v \otimes R^{2}_{r, t}\cdot n\Big) \\
  &=& (\pi_{M, r}\otimes (\pi_{N, r}\otimes \iota_{L_{r}})) (\iota_{M_{s}}\otimes \tau_{r, t}) (\tau_{r, s}\otimes \iota_{N_{t}}) 
      \Big(R^{1}_{r, t}R^{1}_{r, s}\cdot v \otimes R^{2}_{r, s}\cdot m\otimes R^{2}_{r, t}\cdot n\Big) \\
  &=& (\pi_{M, r}\otimes (\pi_{N, r}\otimes \iota_{L_{r}})) (\iota_{M_{s}}\otimes \tau_{r, t}) (\tau_{r, s}\otimes \iota_{N_{t}}) 
      \Big( (R_{r, t})_{1s3} (R_{r, s})_{12t} \cdot (v\otimes m\otimes n)\Big) \\
  &\stackrel{(\ref{2})}{=}& (\pi_{M\otimes N, r}\otimes \iota_{L_{r}}) (\iota_{M_{s}}\otimes \tau_{r, t}) (\tau_{r, s}\otimes \iota_{N_{t}}) 
      \Big( (R_{r, t})_{1s3} (R_{r, s})_{12t} \cdot (v\otimes m\otimes n)\Big).
 \end{eqnarray*}
 Thus if $(\iota \otimes \Delta_{s, t})(R_{r, st}) = (R_{r, t})_{1s3} (R_{r, s})_{12t}$, then 
 $c_{L, M\otimes N} = (\iota_{M}\otimes c_{L, N})(c_{L, M}\otimes \iota_{N})$.
 Conversely, since $\pi_{M\otimes N, r}$, $\tau_{r, t}$ and $\tau_{r, s}$ are bijective,
 $(\iota\otimes\Delta_{s, t})R_{r, st}\cdot (v\otimes m\otimes n)=(R_{r, t})_{1s3} (R_{r, s})_{12t} \cdot (v\otimes m\otimes n)$.
 Let $L=M=N=A$, at this time the module action is the product of $A$, because this product is non-degenerate,
 we can easily get $(\iota\otimes\Delta_{s, t})R_{r, st}=(R_{r, t})_{1s3} (R_{r, s})_{12t}$.
 $\hfill \Box$
 \\
  
 \textbf{Theorem \thesection.4}
 Lat $A=\bigoplus_{p\in G} A_{p}$ be a multiplier Hopf $T$-coalgebra and 
 $R = \{R_{p, q} \in M(A_{p}\otimes A_{q})\}_{p, q\in G}$ be a family of multipliers. 
 Then the monoidal category $({}_{A}\mathcal{M}_{crossed}, \otimes, K, a, l, r)$ is a braided monoidal category with the braiding $c$
 if and only if $A$ is a quasitriangular multiplier Hopf $T$-coalgebra, where $c$ is defined by $R$ as in equation (\ref{2}).

 \emph{Proof} 
 If $c$ is a braiding of the monoidal category, then it follows from Lemmas \thesection.1, \thesection.2 and \thesection.3 that $R$ is a quasitriangular structure.
 Conversely, assume that $R$ is a quasitriangular structure, then by Lemmas \thesection.1, \thesection.2 and \thesection.3, 
 it is enough to show that $c = \{c_{M, N}\}$ is natural.
 
 Let $f: M\rightarrow M'$ and $g: N\rightarrow N'$ be two crossed morphisms of crossed left $A$-$G$-modules, 
 in the following we need to check the braiding preserves the morphisms, i.e., $(g\otimes f)c_{M, N} = C_{M', N'}(f\otimes g)$.
 In fact, for any $m\in M_{s}$ and $n\in N_{t}$ with $\forall s, t\in G$,
 \begin{eqnarray*}
   ((g\otimes f)c_{M, N})_{st}(m\otimes n) 
  &=& (g_{sts^{-1}}\otimes f_{s}) c_{M_{s}, N_{t}} (m\otimes n) \\
  &=& (g_{sts^{-1}}\otimes f_{s}) (\pi_{N, s}\otimes\iota_{M_s}) \tau_{s, t} \big(R_{s,t}\cdot (m\otimes n) \big) \\
  &=& (g_{sts^{-1}}\otimes f_{s})\Big( (\pi_{N, s}\otimes\iota_{M_s})  \big(\tau_{s, t}R_{s,t}\cdot (n\otimes m) \big) \Big) \\
  &=& (g_{sts^{-1}}\otimes f_{s})\Big( (\pi_{s}\otimes\iota)(\tau_{s, t}R_{s,t})\cdot (\pi_{N, s}\otimes\iota_{M_s}) (n\otimes m) \Big) \\
  &=& (\pi_{s}\otimes\iota)(\tau_{s, t}R_{s,t})\cdot \Big( (g_{sts^{-1}}\otimes f_{s})(\pi_{N, s}\otimes\iota_{M_s}) (n\otimes m) \Big) \\
  &=& (\pi_{s}\otimes\iota)(\tau_{s, t}R_{s,t})\cdot \Big( \pi_{N', s}\circ g_{t} (n)\otimes f_{s}(m) \Big),
 \end{eqnarray*}
 and 
 \begin{eqnarray*}
   (C_{M', N'}(f\otimes g))_{st}(m\otimes n) 
  &=&  C_{M'_{s}, N'_{t}} (f_{s}(m)\otimes g_{t}(n)) \\
  &=&  (\pi_{N', s}\otimes\iota_{M'_s}) \tau_{s, t} \Big(R_{s,t}\cdot (f_{s}(m)\otimes g_{t}(n)) \Big) \\
  &=&  (\pi_{N', s}\otimes\iota_{M'_s}) \Big(\tau_{s, t}R_{s,t}\cdot (g_{t}(n)\otimes f_{s}(m)) \Big) \\
  &=&  (\pi_{s}\otimes\iota) (\tau_{s, t}R_{s,t})\cdot \Big( (\pi_{N', s}\otimes\iota_{M'_s})(g_{t}(n)\otimes f_{s}(m)) \Big) \\
  &=& (\pi_{s}\otimes\iota)(\tau_{s, t}R_{s,t})\cdot \Big( \pi_{N', s}\circ g_{t} (n)\otimes f_{s}(m) \Big).
 \end{eqnarray*}
 This completes the proof.
 $\hfill \Box$

\section*{Acknowledgements}

 The work was partially supported by the NNSF of China (No. 11226070), the NJAUF (No. LXY201201019, LXYQ201201103)
 and NSF for Colleges and Universities in Jiangsu Province (No. 11KJB110004).

\end {document}